\documentclass[number,preprint,times]{elsarticle}
\usepackage{amssymb,latexsym,fancyhdr,color}
\usepackage{amsmath, amsthm}

\newtheorem{prop}{\bf Proposition}
\newtheorem{thm}{\bf Theorem}
\newtheorem{corr}{\bf Corollary}

\newcommand{\pf}{\mbox{\bf Proof:}\hspace{3mm}}

\definecolor{titlepagegray}{gray}{0.8}

\begin{document}
\begin{abstract}
In this paper we revisit the integral functional of geometric Brownian motion
$$I_t=  \int_0^t e^{-(\mu s +\sigma W_s)}ds,$$
where $\mu\in\mathbb{R}$, $\sigma>0$, and $(W_s)_{s>0}$ is a standard Brownian motion.

Specifically, we calculate the Laplace transform in $t$ of the cumulative distribution function and of the probability density function of this functional.
\end{abstract}
\makeatletter 
\def\ps@pprintTitle{%
   \let\@oddhead\@empty
   \let\@evenhead\@empty
   \def\@oddfoot{\reset@font\hfil\thepage\hfil}
   \let\@evenfoot\@oddfoot
}
\makeatother
\begin{frontmatter}
\title{Revisiting integral functionals of geometric Brownian motion}

\author[1]{Elena Boguslavskaya}
\ead{elena.boguslavskaya@brunel.ac.uk}
\address[1]{Brunel University, Kingston Ln, London, Uxbridge UB8 3PH, UK}

\author[2]{Lioudmila Vostrikova}
\ead{vostrik@univ-angers.fr}
\address[2]{LAREMA, D\'epartement de
Math\'ematiques, Universit\'e d'Angers,\\ 2, Bd Lavoisier  49045,
\sc Angers Cedex 01, France}


\end{frontmatter}
\begin{section}{Introduction}\label{s1}
\par Assume the canonical filtered probability space $(\Omega, \mathcal{F},\mathbb{F}, P)$ with filtration $\mathbb{F}=(\mathcal{F}_t)_{t>0}$ to satisfy the usual conditions. On this space consider a Brownian motion $X=(X_t)_{t>0}$ with drift $\mu\in \mathbb{R}$ and volatility $\sigma >0$, i.e.
$$X_t= \mu t + \sigma W_t,$$
where $W=(W_t)_{t>0}$ is a standard Brownian motion.
\par We are going to study the integral functional of the corresponding geometrical Brownian motion, namely for $t\geq 0$ we are going to investigate
$$
I_t= \int_ 0^t e^{-X_s} ds= \int_0^t e^{-(\mu s +\sigma W_s)}ds.
\label{I}
$$
The law of the integral functional of geometric Brownian motion of type
$$A_t^{(\mu)}=\int_0^t e^{(2\mu s + 2W_s)}ds$$
was studied by numerous authors.
Alili (1995), Comtet et al.(1998) studied it in the case $\mu=0$. For the case $\mu<0$ it was studied by Comtet and Monthus
(1994,1996). These functionals were also thoroughly studied by Yor (1992a, 1992b,1992c), Schepper et al.(1992), Carmona et al.(1997), Dufresne (2000,2001). In particular, Yor (see 1992a, Proposition 2) states that
$$P\left(A_t^{(\mu)}\in du\,|\, W_t+\mu t=x\right)= \frac{\sqrt{2\pi t}}{u}\exp\left(\frac{x^2}{2t}-\frac{1}{2u}(1+e^{2x})\right)\, \theta_{e^x/u}(t)du$$
where
$$\theta_r(t) = \frac{r}{\sqrt{2\pi^3t}}\,\exp\left(\frac{\pi^2}{2t}\right)\,\int_0^{\infty}\exp\left(-\frac{y^2}{2t}\right)\,\exp(-r\cosh(y))\,\sinh(y)\,\sin\left(\frac{\pi y}{t}\right)\,dy.$$
Dufresne (2000) obtained a series representation for the probability density function of $2A_t^{(\mu)}$ involving generalised Laguerre polynomials and the moments of $2A_t^{(\mu)}$.
Yor (1992c, Theorem 2) showed that
$$2A_{\tau}^{(\mu)}\stackrel{\mathcal{L}}{=} \frac{U}{G},$$
where $\tau$ is independent exponential random variable of the parameter $\lambda$, the variables $U$ and $G$ are independent and distributed as $\mbox{Beta}(1,a_{\mu})$ and $\mbox{Gamma}(b_{\mu},1)$ respectively, with $$a_{\mu}= \frac{\mu+\sqrt{\mu^2+ 2\lambda}}{2}, \; b_{\mu}=a_{\mu}-\mu .$$
Dufresne (2001) showed that the probability density function of $1/\left(2A_t^{(\mu)}\right)$ is given by
$$f_{\mu}(x,t) = e^{-\mu^2t/2}\,p_{\mu}(x,t)$$
with
$$p_{\mu}(x,t) = 2^{-\mu}\,x^{-(\mu+1)/2}\int_{-\infty}^{+\infty} e^{-x\cosh^2(y)}\, q(y,t) \cos\left(\frac{\pi}{2}\left(\frac{y}{t}-\mu\right)\right)\, H_{\mu}\left(\sqrt{x}\sinh(y)\right)dy$$
where $H_{\mu}$ is a Hermite function and
$$q(y,t)= \frac{e^{\pi^2/(8t)-y^2/(2t)}}{\pi\sqrt{2t}}\,\cosh(y).$$
\par In more general setting related to Lévy processes, the following exponential integral functional was intensively studied
\begin{equation}\label{0}
 \int _0^{\infty}\exp(-X_{s-})d\eta _s,
\end{equation}
where $X=(X_t)_{t\geq 0}$ and $\eta= (\eta_t)_{t\geq 0}$ are independent Lévy processes. The conditions for finiteness of integral \eqref{0} were obtained by Erikson and Maller in \cite{EM}. The continuity properties of the law of this integral were studied by Bertoin, Lindner, and Maller in \cite{BLM}.
The equations for the density (under the assumption of existence of smooth densities of these functionals) were provided by Bheme in \cite{Be}, by Bheme and Lindner in \cite{BeL}, and by Kuznetsov, Prado, and Savov in \cite{KPS}.
The  properties of the functional $I_{\tau_q}$ killed at  independent exponential time $\tau_q$   for some parameter $q>0$ were investigated in the papers of Patie and Savov \cite{PS}, and  Prado, Rivero, Van Schaik \cite{PRV}.
\par For fixed time horizon, i.e. for $I_t$,  in the Levy setting for $X$ and $\eta_s=s$, expressions for the Mellin transform, the moments, and the PDE equation for the density  were obtained in Salminen, Vostrikova (2018, 2019) and  Vostrikova (2018).
\par  Such interest to the integral functionals of geometric Brownian motion, and, more generally, to the integral functionals of Levy processes, can be easily explained. These functionals appear in many fields, for example in the study of self-similar Markov processes via Lamperti transform, in the study of diffusions in random environment, in mathematical statistics, in mathematical finance in the evaluation of Asian options, and in the ruin theory.  However,  despite numerous studies, the distributions of  $I_t$ and $I_{\infty}$ are only known for a limited number of cases (cf.\cite{GP}).
\par The main results of this paper are the two explicit expressions (see Theorem 1 and Corollary 2). The first explicit expression is for the Laplace transform of the cumulative distribution function  of the integral functional
of geometric Brownian motion.  The second is for the Laplace transform of the probability density function of the integral functional
of geometric Brownian motion. To our knowledge these results are new.
\par
We proceed in the following way.
Firstly we provide the equation for the probability density of the exponential integral functional of additive processes with fixed time horizon. This result allows us to derive the equation for the probability density function of $I_t$, and to write the equation for its cumulative probability function together with boundary conditions (see Proposition 1). Finally, we derive the equation for the Laplace transform of the complementary cumulative distribution function of $I_t$, relate it to the Kummer equation and solve it explicitly. In Corollary 1 we provide the expressions for the Laplace transform of the cumulative function of $I_t$. In Corollary 2 we provide the expression for the Laplace transform of the probability density function of $I_t$.
\end{section}
\begin{section}{Laplace transform for the cumulative distribution function}
\par Denote by $p_t(x)$, $t>0, x>0$ the probability density function of $I_t$ with respect to Lebesgue measure, and let
$$F(t,y) = P(I_t\leq y) =\int_0^y p_y(x) dx$$
be the cumulative distribution function of $I_t$.
Combining Proposition 2, Proposition 3 and Corollary 2 from \cite{V} we get the following proposition.
\begin{prop}\label{p1}
The law of $I_t$ has a density with respect to Lebesgue measure , and the map $(t,x)\rightarrow p_t(x)$ is of class $C^{\infty}(]0,t],\mathbb{R}^{+,*})$. Moreover, the cumulative distribution function $F(t,y)$ of $I_t$ satisfies the following PDE
\begin{equation}\label{q1}
\frac{\partial}{\partial t}F(t,y)) = \frac{1}{2} \sigma^2 \frac{\partial}{\partial y}(y^2\, \frac{\partial}{\partial y}F(t,y)) - (ay+1)\,  \frac{\partial}{\partial y}F(t,y)
\end{equation}
where $a=\frac{1}{2}\sigma^2-\mu$,\\
with boundary conditions
$F(t,0)=0,\,\,\, \lim_{y\rightarrow +\infty}F(t,y)=1.$
\end{prop}
\par For $t>0$ and $y\geq 0$ define
complementary cumulative distribution function $\bar{F}$
\begin{equation}\label{q2}
\bar{F}(t,y)= 1-F(t,y)
\end{equation}
with  Laplace transform for $\lambda >0$
\begin{equation}\label{q3}
P(y,\lambda)= \int_0^{\infty}e^{-\lambda t}\,\bar{F}(t,y)dt.
\end{equation}
Consider a confluent hypergeometric function of the first kind (Kummer's function) defined as
\begin{equation}\label{M}
M(a,b,z)= \sum_{n=0}^{\infty}\frac{(a)_{n}z^n}{(b)_{n}\,n!}
\end{equation}
where $(a)_{n}$ is a Pochhammer symbol, $(a)_{0}=1, (a)_{n}= a(a+1)(a+2)\cdots (a+n-1)$ and the same for $(b)_{n}$.
\begin{thm} The Laplace transform $P(y,\lambda)$ of $\bar{F}$ satisfies the following differential equation
$$\frac{1}{2}\sigma^2 y^2 P''_{yy}+ (by-1)P'_{y}-\lambda P=0$$
with boundary conditions
$$P(0,\lambda)=\frac{1}{\lambda},_,\,\, \lim_{y\rightarrow +\infty}P(y, \lambda) = 0,$$
or solving it explicitly
\begin{equation}\label{final}
P(y,\lambda) = \frac{1}{\lambda} \left(\frac{2}{y\sigma^2}\right)^{k} \frac{\Gamma\left(1  - \frac{2\mu}{\sigma^2} + k
\right)}
{\Gamma\left(1 - \frac{2\mu}{\sigma^2} + 2k
\right)} M \left(k,1 - \frac{2\mu}{\sigma^2} + 2k, -\frac{2 }{y \sigma^2}\right),
\end{equation}
where $ k = \frac{\mu + \sqrt{\mu^2 + 2 \lambda \sigma^2}}{\sigma^2}$.
\end{thm}
\pf We divide our proof into three parts: firstly we reduce our equation to Kummer's equation and find a general solution, then we adjust this general solution to the boundary  conditions.\\
\it {1) General solution of equation \eqref{q1}. }\rm \\
 From \eqref{q1} and \eqref{q2} we get
\begin{eqnarray}
-\frac{\partial}{\partial t} \bar{F}(t,y) &=& -\frac{1}{2}\sigma^2 \frac{\partial}{\partial y} \left(y^2 \frac{\partial}{\partial y} \bar{F}(t,y)\right) + \left( ay+1\right) \frac{\partial}{\partial y} \bar{F}(t,y),\label{original equation}\\
\bar{F}(t,0)&=&1,\label{boundary at 0} \\
\lim_{y \rightarrow \infty}\bar{F}(t,y)&=&0.\label{boundary at infinity}
\end{eqnarray}
where $a= -\mu + \frac{\sigma^2}{2}$.\\
Expanding the derivative operation and substituting $a= -\mu + \frac{\sigma^2}{2}$ we can rewrite  (\ref{original equation}) as
 \begin{equation}
 \label{PDE}
 \frac{\partial}{\partial t} \bar{F}(t,y) = \frac{1}{2} \sigma^2 y^2 \frac{\partial^2}{\partial y^2} \bar{F}(t,y) + \left( by-1\right) \frac{\partial}{\partial y} \bar{F}(t,y),
 \end{equation}
 where $b=\mu + \frac{\sigma^2}{2}$.
 \par By taking the Laplace transform of (\ref{PDE}) and using \eqref{q3}, we rewrite (\ref{PDE}) as
 \begin{equation}
 \frac{1}{2} \sigma^2 y^2 P_{yy}^{\prime \prime} + \left(by -1\right)P_y^\prime -
 \lambda P = 0 \label{ODE}
 \end{equation}
 From (\ref{boundary at 0}) and from (\ref{boundary at infinity}) we  find  the boundary conditions for $(P(y,\lambda))_{y\geq0,\lambda >0}$:
 \begin{eqnarray}
 \label{P boundary at 0}
 P(0,\lambda) &=& \int_0^\infty e^{-\lambda t} \,\bar{F}(t,0) dt = \int_0^\infty  e^{-\lambda t} dt = \frac{1}{\lambda},\\
 \lim_{y \rightarrow \infty} P(y, \lambda) &=&  \int_0^\infty e^{-\lambda t} \,\left(\lim_{y \rightarrow \infty} \bar{F}(t,y)\right) dt = 0.
 \label{P boundary at infinity}
 \end{eqnarray}
 Next, the equation (\ref{ODE}) can be transformed into
 \begin{equation}
 \frac{1}{2} \sigma^2 \xi u^{\prime \prime}_{\xi \xi} + \left( \xi  + \frac{\sigma^2}{2} -\mu + \sigma^2 k\right) u_\xi^\prime +ku=0.
 \label{eq2}
 \end{equation}
  by setting $y=\xi^{-1}$, $P=\xi^k u$, where $k$ is a root of $\frac{\sigma^2}{2} k^2 - \mu k - \lambda=0$,  i.e.
  \begin{equation}
  \label{k}
  k=\frac{  \mu \pm \sqrt{\mu^2 + 2 \lambda \sigma^2 }}{\sigma^2},
  \end{equation}
  (see eq. 2.1.2.179 from \cite{ZaiPol}).\\
  Equation (\ref{eq2}) is of type 2.1.2.108 in \cite{ZaiPol} and has a solution
  \begin{equation}
  u(\xi) = J\left(k,1- \frac{2\mu}{\sigma^2}+ 2k, -\frac{2  \xi}{\sigma^2}\right),
  \end{equation}
  where $J(a,b;x)$ is any solution of confluent hypergeometric equation
  $$x y^{\prime \prime}_{xx} + (b-x)y^\prime_x -ay=0 $$
  known as Kummer's equation. It is well known there are two fundamental solutions of this equation, namely Kummer's function (confluent hypergeometric function of the first order)
 defined by \eqref{M} and Tricomi's function (confluent hypergeometric function of the second order) defined as
 $$U(a,b,z) =
 \frac{\pi}{\sin(\pi b)}  \left( \frac{M(a,b,z)}{\Gamma(1+a-b) \Gamma(b)} - z^{1-b} \frac{M(1+a-b, 2-b,z)}{\Gamma(a) \Gamma(2-b)}\right).$$
  Therefore, the general solution of the initial problem can be rewritten as
  \begin{equation}
  \label{general solution for P new}
  P(y,\lambda) =c_1 y^{-k} M \left(k,1- \frac{2\mu}{\sigma^2}+ 2k, -\frac{2 }{y \sigma^2}\right) + c_2 y^{-k} U \left(k,1- \frac{2\mu}{\sigma^2}+2k, -\frac{2 }{y \sigma^2}\right),
  \end{equation}
  where $c_1$ and $c_2$ are some real constants.

\it {2) Choice of $k$ and $c_2$ via boundary condition  $\lim_{y \rightarrow \infty} P(y, \lambda) = 0.$}\rm \\
  Note, that there are only two cases for $k$:  $k>0$ if we take the sign $+$ in (\ref{k}), or $k<0$ if we take the sign $-$ in (\ref{k}). Indeed, as $\lambda >0$ we have
  $$
  k = \frac{\mu + \sqrt{\mu^2 + 2 \lambda \sigma^2}}{\sigma^2}  >0,
  $$
  and
  $$
  k = \frac{\mu - \sqrt{\mu^2 + 2 \lambda \sigma^2}}{\sigma^2}  <0.
  $$
In fact only $k>0$ is suitable for our purposes, as both independent solutions explode at $+\infty$ if $k<0$. Moreover, if $k>0$, only the first independent solution is suitable, as the second independent solution also explodes at $+\infty$.
Let us see it in more detail.\\\\
According to  formula 13.5.5, 13.5.10 and 13.5.12 from \cite{AS} for  $a\in \mathbb{R}$ and $b<1$ and $z$ small
\begin{eqnarray*}
M(a,b,z) &=&  1, \mbox{ as   } z \rightarrow 0,\\
U(a,b,z) &=& \frac{\Gamma(1-b)}{\Gamma(1+a-b)}+O\left(|z|^{1-b} \right), \mbox{ for }0< b <1,\\
 &=& \frac{1}{\Gamma(1+a)}  +O\left(|z|\ln(|z|\right), \mbox{ for } b = 0,\\
 &=& \frac{\Gamma(1-b)}{\Gamma(1+a-b)}  +O\left(|z|\right), \mbox{ for } b < 0.
\end{eqnarray*}
Therefore, for $k=\frac{\mu - \sqrt{\mu^2 + 2 \lambda \sigma^2}}{\sigma^2} <0$ we have
$$1- \frac{2\mu}{\sigma^2}+ 2k =1- \frac{2}{\sigma^2}\sqrt{\mu^2 + 2 \lambda \sigma^2}<1,$$
and subsequently
\begin{eqnarray*}
\lim_{y \rightarrow \infty}\left( y^{-k}  M \left(k,1- \frac{2\mu}{\sigma^2}+ 2k, -\frac{2 }{y \sigma^2}\right)\right)&=& \infty\\
\lim_{y \rightarrow \infty}\left( y^{-k}  U \left(k,1- \frac{2\mu}{\sigma^2}+ 2k, -\frac{2 }{y \sigma^2}\right)\right)
&=& \infty.
\end{eqnarray*}
In such a way we know, that for $k<0$ both independent solutions explode, and therefore $c_1$ and $c_2$ should be equal to $0$.
\\\\
 It is easy to check when condition $\lim_{y\rightarrow\infty} P(y,t)=0$ is satisfied for $k>0$.
Indeed, in this case $k = \frac{\mu+ \sqrt{\mu^2+ 2 \lambda \sigma^2}}{\sigma^2}$, and
$$1- \frac{2\mu}{\sigma^2}+ 2k = 1+ \frac{2}{\sigma^2}\sqrt{\mu^2 + 2 \lambda \sigma^2}>1.$$
Thus according to formula 13.5.5 - 13.5.8  in \cite{AS} for  $a\in \mathbb{R}$ and $b>1$ and $z$ small
\begin{eqnarray*}
M(a,b,z) &=&  1, \mbox{ as   } z \rightarrow 0,\\
U(a,b,z) &=& \frac{\Gamma(b-1)}{\Gamma(a)} z^{1-b} +O\left(|z|^{b-2} \right), \mbox{ for } b > 2,\\
&=&  \frac{\Gamma(b-1)}{\Gamma(a)} z^{1-b} +O\left(\ln(|z|) \right), \mbox{ for } b = 2,\\
 &=& \frac{\Gamma(b-1)}{\Gamma(a)} z^{1-b} +O\left(|1|\right), \mbox{ for } 1<  b < 2,
\end{eqnarray*}
we can write
\begin{eqnarray*}
\lim_{y \rightarrow \infty}\left( y^{-k}  M \left(k,1- \frac{2\mu}{\sigma^2}+ 2k, -\frac{2 }{y \sigma^2}\right)\right)=
\lim_{y \rightarrow \infty}\left( y^{-k} M \left(k,1+ \frac{2}{\sigma^2}\sqrt{\mu^2 + 2 \lambda \sigma^2}, -\frac{2 }{y \sigma^2}\right)\right)&=&0,\\
\lim_{y \rightarrow \infty}\left( y^{-k}  U \left(k,1- \frac{2\mu}{\sigma^2}+ 2k, -\frac{2 }{y \sigma^2}\right)\right)=
\lim_{y \rightarrow \infty }\left(y^{-k} U \left(k,1+\frac{2}{\sigma^2}\sqrt{\mu^2 + 2 \lambda \sigma^2}, -\frac{2 }{y \sigma^2}\right)\right)\\
=\lim_{y \rightarrow \infty} \left(y^{-k}\left(\frac{1}{y}\right)^{-\frac{2}{\sigma^2}\sqrt{\mu^2 + 2 \lambda \sigma^2}} \right)=
\lim_{y \rightarrow \infty} \left(y^{\frac{-\mu + \sqrt{\mu^2 + 2 \lambda \sigma^2}}{\sigma^2}} \right)
&=& \infty.
\end{eqnarray*}
In other words only the first independent solution satisfies boundary condition  $\lim_{\lambda \rightarrow \infty} P(y, \lambda) = 0$ when $k>0$, and consequently $c_2$ should be equal to $0$.\\
\it {3)Boundary condition  $P(0, \lambda) = 1/\lambda.$}\rm  \\
According to 13.5.1 in \cite{AS} for large $|z|$ and fixed $a$ and $b$
\begin{eqnarray*}
\frac{M(a,b,z)}{\Gamma(b)}&=& \frac{e^{i\pi a}z^{-a}}{\Gamma(b-a)}\left\{\sum_{n=0}^{R-1}\frac{(a)_n (1+a-b)_n}{n!}(-z)^{-n} + O\left(|z|^{-R} \right)\right\}\\ & &+
\frac{e^z z^{a-b}}{\Gamma(a)}\left\{ \sum_0^{s-1}\frac{(b-a)_n (1-a)_n}{n!}z^{-n}+
O\left(|z|^{-s}\right)\right\}
\end{eqnarray*}
Therefore taking $R=1$ and $s=1$
\begin{eqnarray*}
\lim_{y \rightarrow 0}\left( y^{-k} M \left(k,1- \frac{2\mu}{\sigma^2}+ 2k, -\frac{2 }{y \sigma^2}\right)\right)= \left(\frac{\sigma^2}{2}\right)^k \frac{\Gamma\left(1 - \frac{2\mu}{\sigma^2}+ 2k\right)}{\Gamma\left(1 - \frac{2\mu}{\sigma^2} + k\right)}
\end{eqnarray*}
Finally we get
\begin{equation}
\label{c1c2}
P(0,\lambda) = c_1 \left(\frac{\sigma^2}{2}\right)^k \frac{\Gamma\left(1 - \frac{2\mu}{\sigma^2}+ 2k\right)}{\Gamma\left(1 - \frac{2\mu}{\sigma^2} + k\right)}= \frac{1}{\lambda}.
\end{equation}
and, subsequently,
\begin{equation}
c_1 = \frac{1}{\lambda} \left(\frac{\sigma^2}{2}\right)^{-k} \frac{\Gamma\left(1 - \frac{2\mu}{\sigma^2}+ k\right)}{\Gamma\left(1 - \frac{2\mu}{\sigma^2} + 2k\right)},
\end{equation}
where $k= \frac{\mu + \sqrt{\mu^2 + 2 \lambda \sigma^2}}{\sigma^2}$,
and \eqref{final} is proved.
$\Box$
\corr The Laplace transform $\hat{F}(y,\lambda)$ of the cumulative function $F_t(y)$ of $I_t$ at $\lambda >0$ is given by :
\begin{equation*}\hat{F}(y,\lambda)=
\frac{1}{\lambda}\left\{1- \left( y \frac{\sigma^2}{2}\right)^{-k} \frac{\Gamma\left(1 - \frac{2\mu}{\sigma^2}+ k\right)}{\Gamma\left(1 - \frac{2\mu}{\sigma^2} + 2k\right)} M \left(k,1 - \frac{2\mu}{\sigma^2} + 2k, -\frac{2 }{y \sigma^2}\right)\right\},
\end{equation*}
where $k= \frac{\mu + \sqrt{\mu^2 + 2 \lambda \sigma^2}}{\sigma^2}$.

\pf The result follows directly from the definition of $\bar{F}$ and Theorem 1 since
$\hat{F}(y,\lambda)=\frac{1}{\lambda} - P(y,\lambda).$
$\Box$
\corr The Laplace transform $\hat{p}(y,\lambda)$ of the probability density $p_t(y)$ of $I_t$ at $\lambda >0$ is equal to :\\
\begin{eqnarray*}
\hat{p}(y,\lambda)& = &
\frac{1}{\lambda} \left( y \frac{\sigma^2}{2}\right)^{-k} \frac{\Gamma\left(1 - \frac{2\mu}{\sigma^2}+ k\right)}{\Gamma\left(1 - \frac{2\mu}{\sigma^2} + 2k\right)}
\left\{ \frac{k}{y^{k+1}}\,M \left(k,1 - \frac{2\mu}{\sigma^2} + 2k, -\frac{2 }{y \sigma^2}\right)\right. \\ &    &
\left.
-\frac{2k}{\sigma^2\,y^{k+2}(1 - \frac{2\mu}{\sigma^2} + 2k)}\,
M \left(k+1,2 - \frac{2\mu}{\sigma^2} + 2k, -\frac{2 }{y \sigma^2
}\right)\right\},
\end{eqnarray*}
where $k= \frac{\mu + \sqrt{\mu^2 + 2 \lambda \sigma^2}}{\sigma^2}$.\\
\pf We take the derivative w.r.t. $y$ in the expression of the Laplace transform $\hat{F}(y,\lambda)$ of $F$ and  use 13.4.8 from \cite{AS}
$$\frac{d}{dz} M(a,b,z)= \frac{a}{b}\,M(a+1,b+1,z). \,\,\,\Box $$

\rm Let us denote by $P(y,z), z\in \mathbb{C},$ the extension of the function $P(y,\lambda), \lambda >0,$ constructed in the usual way. Then, since $P(y,z)$ is an analytic function on the half-plan with $Re(z)>0$, the inverse Laplace transform can be calculated by the Bromwich-Mellin formula, namely
$$1-F_t(y) = \frac{1}{2\pi i}\int_{\lambda-i\infty}^{\lambda-i\infty}e^{zt}P(y,z) dz$$
with any $\lambda >0$. The similar formula is valid for the inversion of the Laplace transform  $\hat{p}(y, \lambda)$ of the density $p_t(y)$.
\end{section}
\begin{section}{Acknowledgements}
This research  was partially supported by Defimath project of the Research Federation of "Math\' ematiques des Pays de la Loire" and by PANORisk project "Pays de la Loire" region. We would also like to thank Prof. Michael Bordag from Leipzig University, Germany, for helpful remarks, comments
and numerical calculations related to Bromwich-Mellin formula.
\end{section}

\end{document}